\documentclass{article}

\usepackage{arxiv}

\usepackage[utf8]{inputenc}
\usepackage[T1]{fontenc}
\usepackage{hyperref}
\usepackage{url}
\usepackage{booktabs}

\usepackage{amsmath,amssymb,eucal}

\usepackage{amsfonts}

\title{On the rates of pointwise convergence for Bernstein polynomials}

\newif\ifuniqueAffiliation

\uniqueAffiliationtrue

\ifuniqueAffiliation
\author{Jos\'e A. Adell  \\
	Departamento de M\'etodos Estad\'\i sticos\\
	Universidad de Zaragoza\\
	50009 Zaragoza, Spain \\
	\texttt{adell@unizar.es} \\
	\And
	Daniel C\'ardenas-Morales \\
	Departamento de Matem\'aticas\\
	Universidad de Ja\'en\\
	23071 Ja\'en, Spain \\
	\texttt{cardenas@ujaen.es} \\
    \And
	Antonio J. L\'opez-Moreno \\
	Departamento de Matem\'aticas\\
	Universidad de Ja\'en\\
	23071 Ja\'en, Spain \\
	\texttt{ajlopez@ujaen.es} \\
}
\else

\fi

\renewcommand{\shorttitle}{}

\renewcommand{\headeright}{}

\renewcommand{\undertitle}{}

\hypersetup{pdftitle={},pdfsubject={},pdfauthor={},pdfkeywords={},}

\newtheorem{theorem}{Theorem}

\newtheorem{corollary}{Corollary}

\begin{document}

\maketitle

\begin{abstract}
Let $f$ be a real function defined on the interval $[0,1]$ which is constant on $(a,b)\subset [0,1]$, and let $B_nf$ be its associated $n$th Bernstein polynomial. We prove that, for any $x\in (a,b)$, $|B_nf(x)-f(x)|$ converges to $0$ as $n\rightarrow \infty $ at an exponential rate of decay. Moreover, we show that this property is no longer true at the boundary of $(a,b)$. Finally, an extension to Bernstein-Kantorovich type operators is also provided.
\end{abstract}


\section{Introduction}\label{sec1}
Let $\mathbb{N}$ be the set of positive integers. Denote by $B[0,1]$ the set of all real bounded functions defined on $[0,1]$ endowed with the usual supremum norm $\|\cdot \|$, by $C[0,1]$ the subset of continuous functions, and by $C^k[0,1]$ the subset of $k$-times continuously differentiable functions. Unless otherwise specified, we assume from now on that $n,k\in \mathbb{N}$ and $x\in [0,1]$.

Recall that the $n$th Bernstein polynomial of $f$ is defined as
\begin{equation}\label{1}
B_nf(x)=\sum_{j=0}^nf\left(\frac{j}{n}\right){n \choose j}x^j(1-x)^{n-j},\quad f\in B[0,1].
\end{equation}
The rates of uniform convergence of $B_nf$ towards a function $f\in C[0,1]$, as $n\rightarrow \infty $, are characterized in terms of the so-called Ditzian-Totik second modulus of smoothness of $f$ (see, for instante, \cite{knoopzhou1994}, \cite{knoopzhou1995}, \cite{arXiv2024}, and the references therein).

With respect to the pointwise convergence, many different results have been obtained. Among them, we mention the following (see P\u{a}lt\u{a}nea \cite[pp. 94--96]{paltanea2004})
\begin{equation}\label{2}
\left\vert B_nf(x)-f(x)\right\vert \leq \frac{11}{8}\omega_2\left(f;\frac{\varphi (x)}{\sqrt{n}}\right),\quad f\in C[0,1],
\end{equation}
where $\varphi (x)=\sqrt{x(1-x)}$ and $\omega _2(f;\cdot )$ is the usual second modulus of continuity of $f$ defined as
\[
\omega _2(f;\delta)=\sup\left\{ \left\vert f(x-h)-2f(x)+f(x+h)\right\vert :x\pm h \in [0,1],\ 0\leq h\leq \delta   \right\},\quad \delta \geq 0.
\]
Observe that the order of magnitude of the upper bound in (\ref{2}) cannot be better than $n^{-1}$.

More recently, several authors have obtained a quantitative form of the generalized Voronovskaja's formula (see \cite{gonska}, \cite{tachev}, and \cite{gavreaivan1}, among others). More precisely, let $\mu _{n,j}(x)$ be the $j$th central moment of the Bernstein polynomial, i.e.,
\[
\mu _{n,j}(x)=B_ne_{j,x}(x),\qquad e_{{j,x}}(t)=(t-x)^j,\qquad  t\in [0,1],\ j\in \mathbb{N}\cup \{0\}.
\]
If $f\in C^{2k}[0,1]$, then (cf. \cite[Corollary 3]{jat2018})
\begin{equation}\label{3}
\left\vert B_nf(x)- \sum_{j=0}^{2k}\frac{f^{(j)}(x)}{j!}\mu _{n,j}(x)   \right\vert \leq \frac{1}{6^kn^k}\left(\frac{1}{k!}+\frac{(2r)!}{6^r(k+r)!}\right)\omega \left(f^{(2k)};\frac{1}{\sqrt{n}}\right),
\end{equation}
where $\omega (f^{(2k)};\cdot )$ is the usual first modulus of continuity of $f^{(2k)}$ and
\[
r=\left\lfloor 1+\sqrt{1+\frac{3}{2}k} \right\rfloor ,
\]
$\lfloor \cdot \rfloor$ denoting the integer part. An immediate consequence of (\ref{3}) is that if $f^{(j)}(x)=0$, $j=1, \ldots , 2k$, then the rate of convergence of $|B_nf(x)-f(x)|$ is faster than $n^{-k}$.

On the other hand, let $D_{loc}[0,1]\subset B[0,1]$ be the set of those functions $f$ for which $f^{(i)}(x_0)=0$, $i\in \mathbb{N}\cup \{0\}$, for some $x_0\in (0,1)$. In \cite{antonio}, it was shown that if $f\in D_{loc}[0,1]$, then $B_nf(x_0)=o(n^{-\infty })$, meaning that $B_nf(x_0)=o(n^{-k})$, for all $k\in \mathbb{N}$. This property is actually true for a large class of positive linear operators (cf. \cite{antonio}).

In view of these facts, a natural question is what would be the rate of convergence of $|B_nf(x)-f(x)|$ towards $0$ if $f$ is constant on an interval $(a,b)$ such that $0\leq a<x<b\leq 1$. In the following section, we show that for any $f\in B[0,1]$ (not necessarily continuous) fulfilling this assumption, $|B_nf(x)-f(x)|$ has an exponential rate of decay. This property is based on accurate estimates of the tail probabilities of the binomial distribution. As related remarks, we also show that: (i) the exponents in the rates of convergence cannot be improved, in general; (ii) at the boundary of $(a,b)$, the rate of convergence may be of polynomial order, or even no convergence may occur; and (iii), under the weaker condition that $f\in D_{loc}[0,1]$, an exponential rate of decay at $x_0$ may not be true. Finally, in Section 3, some of the previous results are extended to Bernstein-Kantorovich type operators.

\section{Main results}\label{sec2}
Let $Y_1(x), \ldots , Y_n(x)$ be independent copies of a random variable $Y(x)$ having the Bernoulli distribution with success probability $x$, that is,
\begin{equation}\label{4}
P(Y(x)=1)=1-P(Y(x)=0)=x.
\end{equation}
Denote by $S_n(x)=Y_1(x)+\cdots +Y_n(x)$. Since this random variable has the binomial distribution with parameters $n$ and $x$, we can rewrite (\ref{1}) in probabilistic terms as
\begin{equation}\label{5}
B_nf(x)=\mathbb{E}f\left(\frac{S_n(x)}{n}\right),\quad f\in B[0,1],
\end{equation}
where $\mathbb{E}$ stands for mathematical expectation.

The quantity
\begin{equation}\label{66}
r(x,\theta)=\theta \log \frac{\theta}{x}+(1-\theta)\log \frac{1-\theta}{1-x}, \quad x,\theta \in (0,1),
\end{equation}
is called the Kullback-Leibler divergence between Bernoulli random variables with success probabilities $x$ and $\theta$. Arratia and Gordon \cite{arratiagordon} showed that
\begin{equation}\label{77}
P(S_n(x)\geq bn)\leq e^{-nr(x,b)},\quad 0<x<b<1.
\end{equation}
Since the random variables $S_n(x)$ and $n-S_n(1-x)$ have the same law, we have from (\ref{66}) and (\ref{77})
\begin{equation}\label{88}
P(S_n(x)\leq an)=  P\left(S_n(1-x)\geq n(1-a)\right)   \leq   e^{-nr(1-x,1-a)}=e^{-nr(x,a)},\quad 0<a<x<1.
\end{equation}

Denote by $\lceil y \rceil$ the ceiling of $y\in \mathbb{R}$. Recently, Ferrante \cite{ferrante} has shown the following refinement of (\ref{77}).

\noindent \textbf{Theorem A}
\textit{\ Let $0<x<b<1$. Suppose that
\begin{equation}\label{99}
n\geq 2\qquad \text{and}\qquad 1\leq bn\leq n-1.
\end{equation}
Then,
\begin{equation}\label{100}
P(S_n(x)\geq bn)\leq  \frac{\beta (1-x)}{\beta -x}\frac{1}{\sqrt{2\pi \beta (1-\beta )n}}                  e^{-nr(x,b)},
\end{equation}
where $\beta =\lceil bn\rceil /n$.
}

As in (\ref{88}), an upper bound for the left tail probability $P(S_n(x)\leq an)$, $0<a<x<1$, can be derived, as well.

It is shown in \cite{ferrante} that the upper bound in (\ref{100}) is asymptotically sharp, as $n\rightarrow \infty $. In fact, Theorem A is implicitly an asymptotic result, because, for a fixed $b\in (0,1)$, condition (\ref{99}) implies that
\[
n\geq \max \left(\frac{1}{b},\frac{1}{1-b}\right).
\]

In contraposition to Theorem A, no restrictions on $n$ are required in estimates (\ref{77}) and (\ref{88}). For this reason, we will use such estimates in what follows.

Finally, let $f(x,\theta)$ and $g(x,\theta)$ be two real functions defined on $(0,1)$, for a fixed $x\in (0,1)$. By $f(x,\theta) \sim g(x,\theta)$, $\theta \rightarrow x$, we mean that $f(x,\theta) / g(x,\theta)\rightarrow 1$, as $\theta \rightarrow x$. Differentiating with respect to $\theta $ the function $r(x,\theta)$ defined in (\ref{66}), we obtain $r(x,x)=0$, $r^{(1)}(x,x)=0$, and $r^{(2)}(x,\theta )=1/\varphi ^2(\theta )$. We therefore have
\begin{equation}\label{200}
r(x,\theta ) \sim \frac{1}{2} \frac{(\theta -x)^2}{\varphi ^2 (x)},\quad \text{as} \ \theta \rightarrow x.
\end{equation}
Denote by $1_A$ the indicator function of the set $A$. We state our first main result.

\begin{theorem}\label{teorema2}
Let $f\in B[0,1]$ and $0<a<b<1$. Assume that $f(t)=c$, $t\in (a,b)$, for some real constant $c$. Then, we have for any $x\in (a,b)$
\[
|B_nf(x)-f(x)|\leq \|f-c\|\left(e^{-nr(x,a)}+e^{-nr(x,b)}\right).
\]
\end{theorem}
\
\textit{Proof}. 
Let $x\in (a,b)$. Consider the function $g(t)=f(t)-c$, $0\leq t\leq 1$. Since $g(t)=0$, $t\in (a,b)$, we have from (\ref{5})
\[
|B_nf(x)-f(x)|=|B_ng(x)|=\left\vert \mathbb{E}g\left(\frac{S_n(x)}{n}\right)1_{\{S_n(x)/n\leq a\}}+\mathbb{E}g\left(\frac{S_n(x)}{n}\right)1_{\{S_n(x)/n\geq b\}}\right\vert
\]
\[
\leq \|g\|  \left(  P  \left(  \frac{S_n(x)}{n}\leq a\right)+P\left(\frac{S_n(x)}{n}\geq b \right)\right) \leq \|g\|\left(e^{-n r(x,a)}+e^{-nr(x,b)}\right),
\]
where the last inequality follows from (\ref{77}) and (\ref{88}). The proof is complete.

\hfill $\Box$

In general, there is no hope to obtain a lower inequality for $|B_nf(x)-f(x)|$ in Theorem \ref{teorema2}. Indeed, consider an antisymmetric function $f$ around $1/2$, that is, $f(t)=-f(1-t)$, $t\in [0,1]$. Suppose, in addition, that $f(t)=0$ in $(1/2-\delta ,1/2+\delta )$, for some $0<\delta <1/2$. Since
\[
B_nf\left(\frac{1}{2}\right)=\mathbb{E}f\left(\frac{S_n(\frac{1}{2})}{n}\right)
=-\mathbb{E}f\left(\frac{n-S_n(\frac{1}{2})}{n}\right)=-\mathbb{E}f\left(\frac{S_n(\frac{1}{2})}{n}\right)=-B_nf\left(\frac{1}{2}\right),
\]
we see that $B_nf(1/2)=f(1/2)=0$.

For intervals containing one of the endpoints of $[0,1]$, we give the following results, whose proofs we omit since they follow the lines of that of Theorem \ref{teorema2}.

\begin{corollary}\label{corolario3}
Let $f\in B[0,1]$. Assume that $f(t)=c$, $t\in [0,b)$, $0<b<1$, for some real constant $c$. Then, for any $x\in (0,b)$ we have
\[
|B_nf(x)-f(x)|\leq \|f-c\|e^{-nr(x,b)}.
\]
\end{corollary}

\begin{corollary}\label{corolario4}
Let $f\in B[0,1]$. Assume that $f(t)=c$, $t\in (a,1]$, $0<a<1$, for some real constant $c$. Then, we have for any $x\in (a,1)$
\[
|B_nf(x)-f(x)|\leq \|f-c\|e^{-nr(x,a)}.
\]
\end{corollary}

Concerning the previous results, some remarks are in order. For the sake of concreteness, we focus our attention on Corollary \ref{corolario3}.

\subsection{Sharpness of $r(x,b)$}
Fix $k,m \in \mathbb{N}$, with $1\leq k \leq m-1$, and $m\geq 2$, and let $b=k/m$. Using (\ref{66}) and Stirling's approximation, we have for any $x\in (0,b)$
\[
\left(B_{nm}1_{\{b\}}\right)(x)-1_{\{b\}}(x)=\mathbb{E}1_{\{b\}}\left(\frac{S_{nm}(x)}{nm}\right)=P(S_{nm}(x)=nk)={nm \choose nk}x^{nk}(1-x)^{n(m-k)}
\]
\[
\sim \frac{1}{\sqrt {2nmb(1-b)}}\left(\left(\frac{x}{b}\right)^b\left(\frac{1-x}{1-b}\right)^{1-b}\right)^{nm}=\frac{1}{\sqrt {2nmb(1-b)}}e^{-nmr(x,b)},
\]
as $n\rightarrow \infty $. This shows that the function $r(x,b)$ in the exponent is best possible.

\subsection{Behaviour at the boundary}
Let $(h(n))_{n\geq 1}$ be a sequence of nondecreasing positive real numbers such that $h(n)\rightarrow \infty $ and $h(n)/n\rightarrow 0$, as $n\rightarrow \infty $. Obviously, Corollary \ref{corolario3} is not meaningful for $x=b$, since $r(b,b)=0$, as follows from (\ref{66}). On the other hand, for a fixed $n$, Corollary \ref{corolario3} has interest only if
\[
n r(x,b) \sim \frac{n}{2} \left( \frac{b-x}{\varphi (x)}\right)^2 \geq h(n) \quad \Longleftrightarrow \quad b-x\geq \varphi (x)\sqrt{\frac{2h(n)}{n}},\quad b\rightarrow x,
\]
as follows from (\ref{200}). This suggests that no similar result to Corollary \ref{corolario3} can be given for $x=b$.

In fact, using the convergence of moments in the central limit theorem (see, for instance, Billingsley \cite[p.338]{billingsley}), we have
\begin{equation}\label{20}
\lim _{n\rightarrow \infty}\mathbb{E}\left\vert \frac{S_n(x)-nx}{\sqrt{n}}\right\vert ^{s}=\varphi ^s(x)\mathbb{E}|Z|^s
,\qquad x\in (0,1),\quad s>0,
\end{equation}
where $Z$ is a standard normal random variable. Denote by $y_+=\max (0,y)$, $y\in \mathbb{R}$. Since $S_n(1/2)$ is a symmetric random variable, we have
\begin{equation}\label{20star}
\mathbb{E}\left(\frac{S_n(1/2)}{n}-\frac{1}{2}\right)^s_+=\frac{1}{2}\mathbb{E}\left\vert \frac{S_n(1/2)}{n}-\frac{1}{2}\right\vert ^s,\qquad s>0.
\end{equation}
For any $s>0$, consider the function $f_s\in C[0,1]$ defined as
\[
f_s(t)=\left(t-\frac{1}{2}\right)^s_+,\qquad t\in [0,1].
\]
By (\ref{20}) and (\ref{20star}), we have
\[
B_nf_s\left(\frac{1}{2}\right)-f_s\left(\frac{1}{2}\right)=\frac{1}{2}\mathbb{E}\left\vert \frac{S_n(1/2)}{n}-\frac{1}{2}\right\vert ^s \sim \frac{1}{2^{s+1}}\mathbb{E}|Z|^s\frac{1}{n^{s/2}},\quad n\rightarrow \infty .
\]
In other words, we cannot have an exponential rate of decay in Corollary \ref{corolario3} for $x=b$. Even more, it may happen that $B_nf(b)$ does not converge to $f(b)$. Indeed, let $g\in B[0,1]$ having right and left limits at $t\in (0,1)$, denoted by $g(t^+)$ and $g(t^-)$, respectively. Herzog and Hill \cite{herzoghill} showed that
\[
\lim_{n\rightarrow \infty}B_ng(t)=\frac{1}{2}\left(g(t^+)+g(t^-)\right).
\]
Therefore, for the function $1_{[b,1]}\in B[0,1]$, we have
\[
\lim_{n\rightarrow \infty }\left(B_n1_{[b,1]}\right)(b)=\frac{1}{2}\ne 1=1_{[b,1]}(b).
\]

\subsection{Rates in the set $D_{loc}[0,1]$}
Let $f\in D_{loc}[0,1]$. As said in Section \ref{sec1}, it is known that $B_nf(x_0)=o(n^{-\infty})$. In such circumstances, it is natural to wonder if
\begin{equation}\label{a}
|B_nf(x_0)|\leq Ke^{-C(x_0)n},\quad n\in \mathbb{N},
\end{equation}
for some positive constants $K$ and $C(x_0)$. The following example gives us a negative answer.

Let $0<\alpha <1$. Consider the function
\[
f(t)=\left\{
       \begin{array}{ll}
         \exp (-|t-1/2|^{-\alpha}), & t\in [0,1]\setminus \{1/2\} \\
         0, & t=1/2.
       \end{array}
     \right.
\]
Observe that $f\in D_{loc}[0,1]$ with $f^{(i)}(1/2)=0$, $i\in \mathbb{N}\cup\{0\}$. By Stirling's approximation, we have
\[
B_{2n}f\left(\frac{1}{2}\right)=\mathbb{E}f\left(\frac{S_{2n}(1/2)}{2n}\right)\geq f\left(\frac{n-1}{2n}\right)P(S_{2n}(1/2)=n-1)
\]
\[
=\frac{n}{n+1}f\left(\frac{n-1}{2n}\right)\frac{(2n)!}{(n!)^2}\frac{1}{4^n}\sim \exp \left(-(2n)^{\alpha }\right)\frac{1}{\sqrt{\pi n}},\quad n\rightarrow \infty.
\]
Since $0<\alpha <1$, this shows that inequality (\ref{a}) cannot be true.

\section{Bernstein-Kantorovich type operators}\label{sec4}
Let $k\in \mathbb{N}\cup \{0\}$. Let $W_k$ be a random variable taking values in $[0,k]$ and independent of the random variables $Y_j(x)$, $j=1,\ldots,n$,  considered in Section \ref{sec2}. We define the operator
\begin{equation}\label{21}
L_{n,k}f(x)=\mathbb{E}f\left(\frac{S_{n-k}(x)+W_k}{n}\right),
\end{equation}
for any measurable function $f\in B[0,1]$. If $(V_j)_{j\geq 1}$ is a sequence of independent copies of a random variable $V$ uniformly distributed on $[0,1]$ and $W_k=V_1+\cdots +V_k$, $k\in \mathbb{N}$ ($W_0=0$), then the operator defined in (\ref{21}) is the Bernstein-Kantorovich operator, denoted by $B_{n,k}:=L_{n,k}$. Observe that $B_{n,0}=B_n$. For $k\in \mathbb{N}$, we can write (see Mache \cite{mache}, Acu et al. \cite{acuetall}, and the references therein)
\[
B_{n,k}f(x)=\sum_{j=0}^{n-k}{n-k \choose j}x^j(1-x)^{n-k-j}\int_0^1\cdots \int_0^1f\left(\frac{j+v_1+\cdots +v_k}{n}\right)dv_1\ldots dv_k.
\]
We give the local approximation result for $L_{n,k}$ analogous to Theorem \ref{teorema2}.

\begin{theorem}\label{teorema5}
In the setting of Theorem \ref{teorema2}, assume further that $f$ is measurable and that the interval $I=(na/(n-k),b-k/n)$ is nonempty. For any $x\in I$, we have
\[
|L_{n,k}f(x)-f(x)|\leq \|f-c\|\left(    e^{ -(n-k)r(x,na/(n-k))}  +  e^{-(n-k)r(x,b-k/n)}   \right).
\]
\end{theorem}
\
\textit{Proof}.
Let $x\in I$. Denote
\begin{equation}\label{22}
T_{n,k}(x)=\frac{S_{n-k}(x)+W_k}{n}.
\end{equation}
Let $g$ be as in the proof of Theorem \ref{teorema2}. From (\ref{21}) and (\ref{22}), we have
\[
|L_{n,k}f(x)-f(x)|=|L_{n,k}g(x)|=\left\vert \mathbb{E}g\left(T_{n,k}(x)\right)1_{\{T_{n,k}(x)\leq a\}}  +  \mathbb{E}g\left(T_{n,k}(x)\right)1_{\{T_{n,k}(x)\geq b\}} \right\vert
\]
\begin{equation}\label{23}
\leq \|g\|\left(P(T_{n,k}(x)\leq a)+P(T_{n,k}(x)\geq b)\right).
\end{equation}
From (\ref{22}), we see that
\[
b\leq T_{n,k}(x)\leq \frac{S_{n-k}(x)}{n-k}+\frac{k}{n}.
\]
By virtue of (\ref{77}), this implies that
\begin{equation}\label{24}
P(T_{n,k}(x)\geq b)\leq P\left(\frac{S_{n-k}(x)}{n-k}\geq b-\frac{k}{n}\right)\leq e^{-(n-k)r(x,b-k/n)}.
\end{equation}
Again by (\ref{22}),
\[
a\geq T_{n,k}(x)\geq \frac{n-k}{n}\frac{S_{n-k}(x)}{n-k}.
\]
By (\ref{88}), this entails that
\[
P(T_{n,k}(x)\leq a)\leq P\left(\frac{S_{n-k}(x)}{n-k}\leq \frac{an}{n-k}\right)\leq e^{-(n-k)r(x,an/(n-k))}.
\]
This, (\ref{23}), and (\ref{24}) show the result.

\hfill $\Box$

In the setting of Corollary \ref{corolario3}, the analogous result is the following
\[
|L_{n,k}f(x)-f(x)|\leq \|f-c\|   e^{-(n-k)r(x,b-k/n)},\qquad x\in (0,b-k/n).
\]
A similar statement holds in the setting of Corollary \ref{corolario4}.

\

\textbf{Funding}

The first author is supported by Research Project DGA (E48\_23R). The second and third authors are supported by Junta de Andaluc\'\i a (Research Group FQM-0178).

\

\textbf{Corresponding author}

Correspondence to Daniel Cárdenas-Morales.

\

\textbf{Data availibility}

No data was used for the research described in the article.

\

\textbf{Ethics declarations}

The authors declare that they have no competing interests.

\bibliographystyle{unsrtnat}

\end{document}